\documentclass[reqno,10pt,draft,a4paper]{amsart}
\usepackage{amsmath}
\usepackage{amssymb}
\usepackage{latexsym}

\usepackage[cmtip,matrix,arrow]{xy}
\UseComputerModernTips
\CompileMatrices

\usepackage[dvips]{rotating}

\DeclareMathOperator{\Hom}{Hom}
\DeclareMathOperator{\Ext}{Ext}

\renewcommand{\ge}{\geqslant}
\renewcommand{\le}{\leqslant}

\newcommand{\N}{\mathbb{N}}

\newcommand{\Z}{\mathbb{Z}}

\newcommand{\frob}{^{\mathrm F}}

\newcommand{\rarr}{\rightarrow}

\newcommand{\lrimpl}{\Leftrightarrow}

\DeclareMathOperator{\im}{im}
\newcommand{\Mod}{\mathrm{mod}}

\DeclareMathOperator{\Ind}{Ind}

\newcommand{\ses}{short exact sequence }
\newcommand{\GL}{\mathrm{GL}}
\newcommand{\SL}{\mathrm{SL}}
\newcommand{\St}{\mathrm{St}}

\newcommand{\bj}{\bar{\jmath}}
\newcommand{\bi}{\bar{\imath}}

\newcommand{\sidein}{\begin{turn}{45}$\in$\end{turn}}
\newcommand{\fit}{\phantom{\Big(}}

\newcommand{\bDelta}{\overline{\Delta}}
\DeclareMathOperator{\id}{id}

\begin{document}
\theoremstyle{plain}
\newtheorem{thm}{Theorem}[section]
\newtheorem{propn}[thm]{Proposition}
\newtheorem{cor}[thm]{Corollary}
\newtheorem{clm}[thm]{Claim}
\newtheorem{lem}[thm]{Lemma}
\newtheorem{conj}[thm]{Conjecture}
\theoremstyle{definition}
\newtheorem{defn}[thm]{Definition}
\newtheorem{rem}[thm]{Remark}

\title{Higher Extensions between modules for {$\SL_2$}}
\author{Alison E. Parker}
\address{Department of Mathematics, University of Leicester, Leicester
LE1 7RH, UK}
\email{aep24@mcs.le.ac.uk}



\begin{abstract}
We calculate 
$\mathrm{Ext}^{\bullet}_{\mathrm{SL}_2(k)}(\Delta(\lambda), \Delta(\mu))$, 
$\mathrm{Ext}^{\bullet}_{\mathrm{SL}_2(k)}(L(\lambda), \Delta(\mu))$, 
$\mathrm{Ext}^{\bullet}_{\mathrm{SL}_2(k)}(\Delta(\lambda), L(\mu))$,
and
$\mathrm{Ext}^{\bullet}_{\mathrm{SL}_2(k)}(L(\lambda), L(\mu))$,
where $\Delta(\lambda)$ is the Weyl module of highest weight
$\lambda$, 
$L(\lambda)$ is the simple $\mathrm{SL}_2(k)$-module of highest weight
$\lambda$ and our field $k$ is algebraically closed of positive characteristic.
We also get analogous results for the Dipper-Donkin quantisation.
To do thus we construct the Lyndon-Hochschild-Serre spectral sequence
in a new way, and find a new condition for the $E_2$ page of any
spectral sequence to be the
same as the $E_\infty$ page.
\end{abstract}
\maketitle

\section*{Introduction}
\addtolength{\parskip}{10pt}
Suppose $G$ is a  reductive, semi-simple linear algebraic group over
a field $k$ of positive characteristic $p$.
The category of rational $G$-modules is not very well understood, 
we do not even know the characters of
the simple modules in general. We can however distinguish them by their
highest weight which we use as a labelling. 
The category is not semi-simple and has infinite global dimension.
We do have important objects, the Weyl modules, $\Delta(\lambda)$, and their
duals the induced modules, $\nabla(\lambda)$, with $\lambda$ a dominant
weight.
The simple module,
$L(\lambda)$, is the head of $\Delta(\lambda)$.
Understanding the decomposition numbers of the Weyl modules would
determine the characters of the simple modules.

If we let $G = \GL_n(k)$ and 
truncate to the category to polynomial modules of fixed degree $r$
then we get the module category of a quasi--hereditary algebra --- namely
the Schur algebra. The $\nabla(\lambda)$, $\Delta(\lambda)$, with
$\lambda$ a partition of $r$, are the
costandard and standard objects respectively. 
Homological algebra is an important tool for studying such categories.
In particular, we would like to be able to determine 
$\Ext^*_G(\Delta(\lambda), \Delta(\mu))$,
$\Ext^*_G(\Delta(\lambda), L(\mu))$,
$\Ext^*_G(L(\lambda), \Delta(\mu))$,
and
$\Ext^*_G(L(\lambda), L(\mu))$.
Not many of these groups are known explicitly. A survey of what is
known about $\Ext$ groups for algebraic groups and Hecke algebras may
be found in \cite{coxpar2}.

In this paper we give recursive formulas for these groups for the case $G =
\SL_2(k)$. This, in principle, means that these groups can be calculated
explicitly.
The results for Weyl modules are generalisations of results 
of Cox and Erdmann, \cite{erd1}, \cite{cox1} and
\cite{coxerd} who determined $\Ext^i(\Delta(\lambda), \Delta(\mu))$ for
$\SL_2(k)$ and the quantum group $q$-$\GL_2(k)$ for  $i = 1$ and $2$.
(The $i=0$ case is well known, a proof may be found in 
\cite{coxerd}.) They 
use the Lyndon-Hochschild-Serre spectral and its
associated five term exact sequence
sequence to relate the $G_1$-cohomology to the $G$-cohomology. An
analysis of their results shows that they have proved that the $E_2$
terms are the same as the $E_\infty$ terms if they lie on or below the
diagonal $i +j =2$. 
We use a variant of the Lyndon-Hochschild-Serre
spectral sequence for linear algebraic groups
and show that the $E_2$ pages for the spectral sequences for
the $\Ext$ groups we are interested in
are the same as the $E_\infty$ page. We can then just add up the $q$th
diagonal of the $E_2$ page to get the required $\Ext$ group. 
In particular, we get a nice recursion formulas for some of the $\Ext$
groups and so they may be completely calculated by this method.

The method we use to construct the Lyndon-Hochschild-Serre spectral
sequence is equally valid in the finite group setting where we wish to
calculate cohomology for a finite group 
$G$ using the cohomology of a normal
subgroup $N$ and the quotient group $G/N$, see remark \ref{rem:finite}. We find a new
condition, corollary \ref{cor:E2}, 
on the initial exact couple that implies that the $E_2$ page
is the same as the $E_\infty$ page. Thus section \ref{sect:ss} in the
appropriate finite group setting has applications that may be of
interest to finite group theorists.

We can then use the results of Donkin \cite{donkrat} and 
Cline, Parshall and Scott
\cite{cpstransf} to get some $\Ext$ groups for larger $G$. Let $S$ be the set
of simple roots for $G$, and $\lambda$ and $\mu$ dominant weights for $G$.
Set $m_{\alpha} = \langle \lambda, \alpha \rangle$ and 
$n_{\alpha} = \langle \mu, \alpha \rangle$ for $\alpha \in S$.
If $\lambda -\mu =m \beta$ for a $\beta \in S$ and $m \in \Z$ then
$$\Ext^q_G(\Delta(\lambda), L(\mu)) \cong
\Ext^q_{\SL_2}(\Delta(2m_\beta),L(2n_\beta))$$
using \cite[corollary 10]{cpstransf} and
$$\Ext^q_G(\Delta(\lambda), \Delta(\mu)) \cong
\Ext^q_{\SL_2}(\Delta(2m_\beta),\Delta(2n_\beta))$$
using \cite[section 4]{erd1}. 
For $\GL_n(k)$ the condition on $\lambda$ and $\mu$ is equivalent to
saying that they only differ in two consecutive rows.

We finish off by showing that the results for $\SL_2(k)$
may be easily generalised to
the quantum group $q$--$\GL_2(k)$ of Dipper and Donkin \cite{dipdonk},
thus analogous results hold there as well.

\section{Notation}
We first review most of the notation that 
we will be using. The reader is referred to
\cite{humph} and \cite{springer} for further information. 
This material is also in~\cite{jantz2} where is it presented in 
the form of group schemes.

Throughout this paper $k$ will be an algebraically closed field of
characteristic $p$.
Let $G$ be a linear algebraic group which is connected and
reductive, and let  $\mathrm{F}:G \rarr G$ its  corresponding
Frobenius morphism. We let $G_1$ be the first Frobenius
kernel.
We fix  a maximal torus $T$ of $G$ of dimension $n$, 
the rank of $G$.
We also fix $B$, a Borel subgroup of $G$ with $B \supseteq T$ 
and let $W_p$ be the affine Weyl group of $G$.

We will write $\Mod(G)$ for the category of finite
dimensional rational $G$-modules. Most $G$-modules considered in this
paper will belong to this category. 
Let $X(T)=X$ be the weight lattice for $G$. 
We take $R$ to be the roots of $G$ and $S$ the simple roots.
Let $h$ be the maximum of the Coxeter numbers for the connected
components of $R$.

We have a partial order on $X$ defined by 
$\mu \le \lambda \lrimpl \lambda -\mu \in \N S$.
We let $X^+$ be the set of dominant weights.
Take $\lambda \in X^+$ and let $k_\lambda$ be the one-dimensional module
for $B$ which has weight $\lambda$. We define the \emph{induced
module}, $\nabla(\lambda)= \Ind_B^G(k_\lambda)$. 
This module has formal character given by Weyl's character formula and has
simple socle $L(\lambda)$, 
the irreducible $G$-module of highest weight
$\lambda$.  Any finite dimensional, rational irreducible $G$-module
is isomorphic to $L(\lambda)$ for a unique $\lambda \in X^+$.

Since $G$ is split, connected and reductive we have an
antiautomorphism, $\tau$, which acts as the identity on $T$
(\cite{jantz2}, II, corollary 1.16). From this morphism we 
may define $^\circ$, 
a contravariant dual. It does not change a module's character, hence
it fixes the irreducible modules. 
We define the Weyl module, to be 
$\Delta(\lambda)=\nabla(\lambda)^\circ$.
Thus $\Delta(\lambda)$ has simple head $L(\lambda)$.
We say a module is a \emph{tilting module} if it has a filtration by
both Weyl modules and induced modules. For each $\lambda \in X^+$
there is a unique indecomposable tilting module $T(\lambda)$ with
highest weight $\lambda$.

We say that $\lambda$ and $\mu$ are \emph{linked} if  they belong to the 
same $W_p$ orbit on $X$ (under the dot action).
If two irreducible modules $L(\lambda)$ and $L(\mu)$ are in the same $G$
block then $\lambda$ and $\mu$ are linked. 
We will identify a block with the set of dominant weights which label
the simples in that block.

We let $X_1$ be the $p$-restricted weights. Then the $G$-modules
$L(\lambda)$ with $\lambda \in X_1$ form a complete set of simples
upon restriction for $G_1$.

If $G=\SL_2(k)$ then $X^+$ may be identified with the natural
numbers. Two regular weights $\lambda=pa+i$ and $\mu=pb+j$ with $a$, $b \in \N$
and $0 \le i \le p-2$ and $0 \le j \le p-2$ will only be in the same
block if $a-b$ is even and $i=j$ or if $a-b$ is odd and $j=p-2-i$.
The only non-regular weights are the Steinberg weights, two such
weights $pa+p-1$ and $pb+p-1$ will only be in the same block if $a$
and $b$ are in the same block.

If $H=G$ or $G_1$ or $G/{G_1}$ then 
the category of rational $H$-modules has enough injectives and so we may
define $\Ext_H^*(-,-)$ as usual by using injective
resolutions (see \cite{benson}, section 2.4 and 2.5).
We also define the cohomology groups $H^*(H,-)\cong\Ext^*_H(k,-)$,
which are the right derived functors of the fixed point functor $-^H$.

We may form the Frobenius twist of a module $V$, by composing the given
action with the Frobenius morphism. This new module is denoted
$V\frob$ and it is trivial as $G_1$-module. Conversely any
$G$-module $W$ which is trivial as a $G_1$-module, is of the form 
$V\frob$ for some $G$-module $V$ which is unique up to isomorphism. 
If $W$ and $V$ are $G$-modules then $\Ext^i_{G_1}(W,V)$ has a natural
structure as a $G$-module.
Moreover when $W$ and $V$ are finite dimensional we have,
\begin{equation}\label{g1triv}
\Ext^i_{G_1}(W,V\otimes Y\frob)\cong \Ext^i_{G_1}(W,V) \otimes Y\frob
\end{equation}
as $G$-modules.
If $H=G$ or $G_1$ or $G/{G_1}$ and $V$, $W$ are $H$-modules then we have
$$
\Ext^i_H(W,V)\cong \Ext^i_H(V^*,W^*)\cong \Ext^i_H(k,W^*\otimes
V) \cong
H^i(H,W^*\otimes V)
$$
where $^*$ is the ordinary dual. 
We also note that $(V\frob)^{G/G_1}\cong
V^G$, and $H^i(G/G_1,V\frob)\cong H^i(G,V)$.

\section{Spectral sequences}\label{sect:ss}
We will use a variant of the Lyndon-Hochschild-Serre
spectral sequence to show that the $E_2$ page
is the same as the $E_\infty$ page for the $\Ext$ groups we are
interested in.
Our primary reference for this section is
\cite[chapter 3]{benson2}.

Usually we would construct the Lyndon-Hochschild-Serre spectral
sequence
by taking $G$ and $G/G_1$ injective resolutions as in
\cite[I, section 6.6]{jantz2}. We wish to
analyse the $E_0$ and $E_1$ page so we construct it in a different
way.

First take a $G$ (co)resolution of a finite dimensional $G$-module 
$W$ which is injective as a
$G_1$ resolution. This is certainly possible as the injectives for $G$
are also injective for $G_1$ upon restriction. 
Usually we will want much smaller,
finite dimensional modules in our resolution. Thus the modules will
not be injective as $G$-modules. This will sometimes be
possible as the indecomposable 
$G_1$-modules $Q(i)$ have a $G$ structure (proved in general for $p
\ge 2(h-1)$ \cite{jandar}, and is true for all $p$ for $\SL_2$ and $\SL_3$).
In all the examples we take for $\SL_2$ we will construct this
resolution explicitly.
So we have a $G_1$-injective resolution with $I_i\in\Mod(G)$
$$
0 \to W \to  
I_0 \to I_1 \to \cdots \to I_{m-1} \to I_m \to \cdots.
$$

Now take a (minimal) injective $G$ resolution of $k$.
$$ 0 \to k \to J_0 \to J_1 \to \cdots $$
These injective modules will have infinite dimension. 
We then take the Frobenius twist of the modules $J_m$. We end up with a
$G/G_1$ resolution of $k$ which is injective as $G/G_1$-modules.

We take a finite dimensional $G$-module $V$ and 
form a double complex in a  similar way to \cite[section 3.5]{benson2}
$$E_0^{mn} = \Hom_{G/G_1}(k, \Hom_{G_1}(V,I_n)\otimes J_m\frob)
$$
where the horizontal maps are the differentials induced by the
differentials in the $G/G_1$ resolution of $k$ and the
vertical maps are $(-1)^m$ times the differentials induced by the
$G_1$ injective resolution of $V$.

Using \cite[theorem 3.4.2]{benson2} we now have a spectral sequence with 
\begin{align*}
E_1^{mn}&= \Hom_{G/G_1}(k, \Ext^n_{G_1}(V,W)\otimes J_m\frob)\\
E_2^{mn}&= H^m({G/G_1}, \Ext^n_{G_1}(V,W))\\
\end{align*}

We have 
\begin{align*}
\Hom_{G/G_1}(k, \Hom_{G_1}(V, I_n)\otimes J_m\frob)
&\cong
((V^* \otimes I_n)^{G_1} \otimes J_m\frob)^{G/G_1}\\
&\cong
((V^* \otimes I_n \otimes J_m\frob)^{G_1})^{G/G_1}\\
&\cong
(V^* \otimes I_n \otimes J_m\frob)^{G}\\
&\cong
\Hom_{G}(V, I_n \otimes J_m\frob)\\
\end{align*}
as $J_m\frob$ is trivial as a $G_1$-module and $I_n$ has a
$G$-structure.

Now 
$$
\Ext^q_G(V, I_n \otimes J_m\frob)
\cong 
\Ext^q_{G/G_1}(k, \Hom_{G_1}(V, I_n)\otimes J_m\frob)
$$
as $V^* \otimes I_n \otimes J_m\frob$ is injective as a
$G_1$-module and the usual Lyndon-Hochschild-Serre spectral sequence
collapses to a line.
But now this $\Ext$ group is zero for $q \ge 1$, as $J_m \frob$ is
injective as a $G/{G_1}$-module.
Thus the total homology of the double complex $E_0$
is $\Ext^*_G(V,W)$.

Hence we have constructed a spectral sequence 
which converges to $\Ext_G^q(V,W)$ and whose
$E_2$ page has the same terms as that as the usual
Lyndon-Hochschild-Serre sequence.
(It will also have the same differentials - but we will not need
this.)

Our results for $\SL_2(k)$ say that the $E_2$ page is often the same
as the $E_\infty$ page. We now present some lemmas which give
conditions for the $k_2$ associated to the 2nd derived couple to be
zero.
For what follows we 
refer the reader to \cite[section 3.4]{benson2}, where the $E_1$
and $E_2$ pages as well as the $D_0$ and $D_1$ are explicitly
constructed.

We have the double complex
$$
\xymatrix@R=15pt@C=15pt{
{\vdots}         & {\vdots} & {\vdots}
\\
E_0^{02}\ar@{->}^{d_1}[r]\ar@{->}_{d_0}[u] 
&E_0^{12}\ar@{->}^{d_1}[r]\ar@{->}_{d_0}[u] 
&E_0^{22}\ar@{->}^{d_1}[r]\ar@{->}_{d_0}[u] 
& {\cdots}
\\
E_0^{01}\ar@{->}^{d_1}[r]\ar@{->}_{d_0}[u] 
&E_0^{11}\ar@{->}^{d_1}[r]\ar@{->}_{d_0}[u] 
&E_0^{21}\ar@{->}^{d_1}[r]\ar@{->}_{d_0}[u] 
& {\cdots}
\\
E_0^{00} \ar@{->}^{d_1}[r]\ar@{->}_{d_0}[u] 
&E_0^{10}\ar@{->}^{d_1}[r]\ar@{->}_{d_0}[u] 
&E_0^{20}\ar@{->}^{d_1}[r]\ar@{->}_{d_0}[u] 
& {\cdots}
}
$$
with $d_1^2 =0$, $d_0^2=0$ and $d_1d_0 +d_0d_1 =0$.

We have 
$$D_0^{mn}= \bigoplus_{m+n=e+f, \ e \ge m}E^{ef}_0$$
$$E_1^{mn}= H(E^{mn}_0, d_0)$$
$$D_1^{mn}= H(E^{mn}_0 \oplus E_0^{m+1,n-1}\oplus\cdots, d_0+d_1)$$
We will use square brackets to denote the class of $(x_1, x_2,
\ldots)$ in the homology group.

The first derived couple has long exact sequence 
$$ 
\cdots \to E_1^{m,n-1} \xrightarrow{k_1^{m,n-1}}
D_1^{m+1,n-1} \xrightarrow{i_1^{m+1,n-1}} 
D_1^{m,n} \xrightarrow{j_1^{m,n}}
  E_1^{m,n} \xrightarrow{k_1^{m,n}} 
\cdots
$$
where $k_1[x]=[(d_1x,0,\ldots)]$, 
which is induced by
taking the homology of the short exact sequence
$$0\to D_0^{m+1,n-1} \to D_0^{m,n} \to E_0^{m,n} \to 0$$

We define the higher derived couples by taking the derived couple of
the previous one.
We have an exact diagram of doubly graded $k$-modules
$$
\xymatrix@R=20pt@C=10pt{
D_l \ar@{->}[rr]^{i_l}& &D_l \ar@{->}[dl]^{j_l}\\
& E_l \ar@{->}[ul]^{k_l} &
}
$$
The derived couple (for $l \ge 1$) is defined by
\begin{align*}
D_{l+1}^{mn}&= \im{i_l^{m+1, n-1}}\subseteq D_l^{mn}\fit \\
E_{l+1}^{mn}&= H(E^{mn}_l, d_l)\fit \\
i_{l+1}^{mn}&=i_l^{mn}\Bigl\vert_{D_{l+1}} \\
j_{l+1}^{mn}(i_l^{m+1,n-1}(x))&= j_l^{m+1,n-1}(x) + \im(d_l)\fit  \\
k_{l+1}^{mn}(z +\im(d_l))&= k_l^{mn}(z)\fit \\
d_{l+1} &= j_{l+1} \circ k_{l+1}\fit 
\end{align*}

Note there is a slight abuse of notation here --- we have two maps
called $d_1$. We will only use $d_1$ to denote the horizontal
differentials on the $E_0$ page and use $j_1k_1$ to refer to the
differential of the $E_1$ page. The connection between the two maps is
that $j_1k_1[x] = [d_1x]$.
Since $E_2^{m,n}$ is the kernel of
$j_1k_1$ over its image, we thus have that 
the kernel of $j_1k_1 : E_1^{m,n}\to E_1^{m+1,n}$
is $\{ [x] \mid d_1x=d_0z \mbox{ for some }z \in E_1^{m+1,n-1} \}$.

Since we have a first quadrant spectral sequence the $E_l^{m,n}$ and
$D_l^{m,n}$ must eventually stabilise,
the $E_\infty$ and $D_\infty$ are then this stable
value.

\begin{lem}\label{lem:dzero}
Suppose $d_0^{m+1,n-1}$ is zero, then $k_2^{m,n}$ is zero.
\end{lem}
\begin{proof}
Now since $d_0^{m+1,n-1}$ are zero, we must have that the kernel of
$j_1k_1:E_1^{m,n} \to E_1^{m+1,n-1}$ is the set 
$$\{[x] \mid x \in \ker(d_1:E_0^{m,n} \to
  E_0^{m+1,n}) \}.$$
Thus we have 
$$k_2^{m,n}([x] +\im{j_1k_1})
= k_1^{m,n}([x])=[(d_1x,0,0,\ldots)] = 0$$ 
as $[x] \in \ker{j_1k_1}$ and so $d_1x =0$.
\end{proof}

\begin{lem}\label{lem:dinj}
Suppose $d_0^{m+2,n-1}$ is injective, then $k_2^{m,n}$ is zero.
\end{lem}
\begin{proof}
The kernel of $j_1k_1 : E_1^{m,n}\to E_1^{m+1,n}$
is $\{ [x] \mid d_1(x)=d_0(z) \mbox{ for some }z \in E_1^{m+1,n-1}
\}$.
Thus we have 
$$k_2^{m,n}([x] +\im{j_1k_1})
= k_1^{m,n}([x])=[(d_1x,0,0,\ldots)]$$ 
Now since $[x] \in \ker{j_1k_1}$, $d_1x = d_0 z$ for some $z \in
E_1^{m+1,n-1}$. Thus
$$
[(d_1x,0,0,\ldots)] 
= [(d_0z,0,0,\ldots)] 
= [(d_0z,0,0,\ldots) -(d_1+d_0)(z,0,0,\ldots)]
= [(0,-d_1z,0,\ldots)]. 
$$ 
Now since $d_0:E_0^{m+2,n-1} \to E_0^{m+2,n}$ is an embedding,
we have $d_1z = 0$ if and only if $d_0d_1z =0$ which is zero if and
only if $-d_1d_0 z=0$. But $d_1d_0z=d_1 d_1 x=0$. Thus $d_1z=0$ and 
$k_2^{m,n}$ is zero.   

This is illustrated in the following diagram.
$$
\newdir{(-}{{\begin{turn}{45}{$\in$}\end{turn}}}
\xymatrix@!0{
&&{\vdots} &&&{\vdots} &&&{\vdots} &
\\
&&&&&&&&&&
\\
\cdots \ar@{->}[rr]& &{E_0^{m,n}}\ar@{->}^-{d_1}[rrr]\ar@{->}_{d_0}[uu] 
&&&E_0^{m+1,n}\ar@{->}^-{d_1}[rrr]\ar@{->}_{d_0}[uu] 
&&&{E_0^{m+2,n}}\ar@{->}^-{d_1}[rr]\ar@{->}_{d_0}[uu] 
& &{\cdots}
\\
&x \ar@{|->}'[r][rrr]\ar@{}[ur]|{\!\! \sidein} 
& &&d_1x=d_0z \ar@{|->}[rrr] \ar@{}[ur]|{\!\! \sidein} 
&& & 0 \ar@{}[ur]|{\!\!\!\! \sidein} & &
\\
&&&&&&&&&&
\\
\cdots \ar@{->}[rr] &&E_0^{m,n-1} \ar@{->}^-{d_1}[rrr]\ar@{->}_{d_0}[uuu] 
&&&{E_0^{m+1,n-1}}\ar@{->}^-{d_1}[rrr]\ar@{->}_{d_0}[uuu] 
&&&{E_0^{m+2,n-1}}\ar@{->}^-{d_1}[rr]\ar@{->}_{d_0}[uuu] 
& &{\cdots}
\\
& & &&z \ar@{|->}'[u][uuu] \ar@{|->}'[r][rrr] \ar@{}[ur]|{\!\!\!\!\!\!
  \sidein} 
& &&d_1z  \ar@{|->}'[u][uuu] \ar@{}[ur]|{\!\!\!\!\!\!\!\! \!\! \sidein}  & &
\\
&&{\vdots} \ar@{->}_{d_0}[uu] 
&&& {\vdots}\ar@{->}[uu]  
&&&{\vdots}\ar@{->}_{d_0}[uu] 
}
$$
\end{proof}

\begin{cor}\label{cor:E2}
Suppose the $d_0$ are all either injections or all zero. Then the $E_2$
page is the same as the $E_\infty$ page.
\end{cor}
\begin{proof}
The conditions on the $d_0$ imply that 
all the $k_2$'s in the second derived couple are zero. This in
particular implies that all subsequent differentials $d_f$ are zero
for $f \ge 2$. Since the $E_{i+1}$ is the homology of $E_i$ with
respect to $d_i$ this means that the homology stabilises with the $E_2$
page. Thus the $E_2$ page is the same as the $E_\infty$ page.
\end{proof}

\begin{rem}
The first condition of this corollary, that all the $d_0$ are zero
implies that the $E_2$ page is the same as the $E_\infty$ page is used
implicitly in Evens' proof of a theorem of Nakaoka, 
\cite[theorem 5.3.1]{evens}.
I am grateful to David Benson for this remark.

This corollary and the preceeding lemmas hold for \emph{any} spectral
sequence that is constructed using the method of derived couples. 
\end{rem}

\begin{rem}\label{rem:finite}
The method used to construct the Lyndon-Hochschild-Serre spectral
sequence in this section is completely general.
Thus if we replace the linear algebraic group $G$ by a finite group
$G$
and $G_1$ by a normal subgroup $N$ we get a new construction of the 
Lyndon-Hochschild-Serre spectral sequence for finite groups.
\end{rem}

\section{The $\SL_2$ case}\label{sect:sl2g1}
We now restrict our attention to $G = \SL_2$.
We list some general facts about modules for $\SL_2$.
Here $X^+$ may be identified with $\N$, and $X_1$ is the set $\{i \mid
0 \le i \le p-1, i \in \N \}$.

Take $i$, $j \in \N$ with $0\le i \le p-2$ and
$0\le j \le p-2$.
We define $\bi= p-2-i$ and similarly
$\bj= p-2-j$. Note that if $p \ge 3$ means that $i
\ne \bi$. 
The simples $L(i)$ and $L(\bi)$ are in the same $G_1$-block and are
the only simples in this block.
If $p=2$ then $L(i)=L(\bi)$ and there is only one simple in this
$G_1$-block.
We will often need to argue separately for the $p=2$ case for this
reason, as the $G_1$ cohomology is more complicated.
Although the results for $p=2$ can be thought of as the results for
$p\ge 3$ but with the $i=j$ and $i=\bj$ cases added together.

For this section we will assume that $a$, $b$, $i$ and $j \in\N$ and 
$0\le i \le p-2$ and $0\le j\le p-2$.
Thus if $p=2$ then $i=j=0=\bi=\bj$.

We know that the tilting module $T(pa+i)$ for $\SL_2(k)$ is isomorphic
to $T(a-1)\frob \otimes T(p+i)$
for $a \ge 1$.
We also know that  $T(i) \cong \Delta(i) \cong \nabla(i)$. 
We have that 
$T(p+i)=Q(\bi)$, the $G_1$-injective hull of $L(\bi)$  
as a $G$-module.
We have the short exact sequences 
\begin{equation}\label{xanth2}
0 \rarr \Delta(pa+i) 
\rarr \Delta(a-1)\frob \otimes T(p+i)
\rarr \Delta(p(a-1)+\bi) \rarr 0
\end{equation}
and
\begin{equation}\label{xanth1}
0 \rarr \Delta(a-1)\frob \otimes L(\bi)
\rarr \Delta(pa+i)
\rarr \Delta(a)\frob \otimes L(i)\rarr 0
\end{equation}
which appear in the thesis \cite{xanth}.
The latter non-split sequence is a $G_1$ socle series for $\Delta(pa+i)$.
That is $\Delta(a) \frob \otimes L(i)$ is the $G_1$ head of
$\Delta(pa+i)$ and $\Delta(a-1)\frob \otimes L(\bi)$ is the $G_1$
socle, if $a\le1$. 

We also have:
$$\Delta(pa+p-1) \cong \Delta(a)\frob \otimes \St$$
$$T(pa+p-1) \cong T(a)\frob \otimes \St$$
where $\St= L(p-1)$ is the Steinberg module,
and natural isomorphisms
$$\Ext_G^m(M,N) 
\cong \Ext_G^m(M\frob \otimes \St, N\frob \otimes \St)$$
for all $m\in \N$ and $M$, $N$ $G$-modules. 
  
The aim is to calculate $\Ext^m_G(\Delta(pb+j),\Delta(pa+i))$,
$\Ext^m_G(L(pb+j), \Delta(pa+i))$,
$\Ext^m_G(\Delta(pb+j), L(pa+i))$,
$\Ext^m_G(L(pb+j), L(pa+i))$,
and $\Ext^m_G(T(pb+j), \Delta(pa+i))$.

First note that for $N$ and $M$ $G$-modules and $q \ge 0$ we have
$$
\Ext^q_{G_1}(N\frob \otimes L(i), M\frob \otimes Q(i))
\cong \left\{ \begin{array}{ll}
           (N\frob)^* \otimes M\frob   &\mbox{if $q=0$}\\ 
           0                           &\mbox{otherwise.} 
              \end{array} \right.
$$
and if $i \ne \bi$ then
$$
\Ext^q_{G_1}(N\frob \otimes L(i), M\frob \otimes Q(\bi))
\cong 0.
$$
as $Q(i)$ is injective as a $G_1$-module.
This will be used extensively in this paper without further comment.

We will use the following lemmas in our calculations in the following
sections. The quantum version of the next two lemmas appears in
\cite[proposition 2.1]{coxerd}.
\begin{lem}\label{lem:Mfrob}
Let $M$ be a $G$-module. Then if $p\ge 3$ we have
$$
\Hom_{G_1}(\Delta(pb+j), M\frob \otimes Q(i))
\cong \left\{ \begin{array}{ll}
           \nabla(b)\frob \otimes M\frob  &\mbox{if $i=j$}\\ 
	   \nabla(b-1)\frob\otimes M\frob &\mbox{if $i=\bj$}\\
           0                              &\mbox{otherwise.} 
              \end{array} \right.
$$
If $p=2$ then
$$
\Hom_{G_1}(\Delta(2b), M\frob \otimes Q(0))
\cong \nabla(b)\frob \otimes M\frob 
	   \oplus \nabla(b-1)\frob\otimes M\frob
$$
\end{lem}
\begin{proof}
We need only prove this for $M=k$ as the $M\frob$ comes out using
\eqref{g1triv}.
Applying the exact functor $\Hom_{G_1}(-, Q(i))$ to the short exact 
sequence \eqref{xanth1} 
gives us 
\begin{multline*}
0 \to
\Hom_{G_1}(\Delta(b)\frob \otimes L(j),Q(i))
\to
\Hom_{G_1}(\Delta(pb+j),Q(i))\\
\to
\Hom_{G_1}(\Delta(b-1)\frob\otimes L(\bj),Q(i))\to 0
\end{multline*}
The first $\Hom$ group is $0$ if $i\ne j$ and $\nabla(b)\frob$ if
$i=j$.
The last $\Hom$ group is $\nabla(b-1)\frob$ if $i=\bj$ and $0$ if
$i\ne \bj$.

If $p=2$ and $i=j$ then both the first and third $\Hom$ groups are
non-zero.
The middle $\Hom$ group has the structure of a $G/G_1$-module.
We have
$$
\Ext^1_{G/G_1}(\nabla(b-1)\frob,\nabla(b)\frob)
\cong 
\Ext^1_{G}(\nabla(b-1),\nabla(b))
\cong 0.
$$
So the middle $\Hom$ group is a direct sum as claimed.
\end{proof}

\begin{lem}\label{lem:MT}
Let $M$ be a $G$-module and  $q \in \N$. If $p\ge3 $ then
$$
\Ext^q_G(\Delta(pb+j), M \frob \otimes Q(i))
\cong \left\{ \begin{array}{ll}
          \Ext^q_G(\Delta(b), M) &\mbox{if $i=j$}\\ 
	  \Ext^q_G(\Delta(b-1), M) &\mbox{if $i=\bj$}\\
          0           &\mbox{otherwise.} 
              \end{array} \right.
$$
If $p=2$ then
$$
\Ext^q_G(\Delta(2b), M \frob \otimes Q(0))
\cong \Ext^q_G(\Delta(b), M) 
	\oplus \Ext^q_G(\Delta(b-1), M)
$$
\end{lem}
\begin{proof}
$$\Ext^m_G(\Delta(pb+j), M \frob \otimes Q(i))
\cong H^m(G/G_1, \Hom_{G_1}(\Delta(pb+j), M \frob \otimes Q(i)))
$$
as the module $\nabla(pb+j)\otimes T(p+i) \otimes M \frob$ is
injective as a $G_1$-module so the usual Lyndon-Hochschild-Serre 
spectral sequence for this module collapses to a line.
We then use the previous lemma.
\end{proof}

The following observation and proof is due to Stephen Donkin.
\begin{lem}\label{lem:good}
The module $\Delta(a)\otimes \nabla(b)$ has a good filtration if $b
\ge a-1$ and it has a Weyl filtration if $b \le a+1$.
\end{lem}
\begin{proof}
Suppose that $b \ge a-1$. The weights of $\Delta(a)$ are
$a$, $a-2$, $a-4$, $\ldots$ , $-a+2$ and $-a$. So the weights of
the $B$-module $\Delta(a)\otimes k_b$ are $a+b$, $a+b-2$, $\ldots$,
$b-a+2$ and $b-a$. Since $b-a \ge -1$ when we induce this module
up to $G$ we will get $\Delta(a)\otimes \nabla(b)$ and it will have a
good filtration (starting at the top) by $\nabla(a+b)$,$ \nabla(a+b-2)$,
$\ldots$, $\nabla(b-a+2)$ and $\nabla(b-a)$, where if $b-a=-1$ then
the last module is taken to be the zero module.
The case $b \le a+1$ is proved by taking duals.
\end{proof}

\begin{lem}\label{lem:dextone}
\begin{align*}
\Ext^1_{G_1}(\Delta(pb+i), \Delta(pa+j)) 
           &{\cong \left\{\begin{array}{ll}
                \Ext_{G_1}^1(\Delta(p(b-a)+i), \Delta(j)) 
                               &\mbox{if $b\ge a$}\\
                \Ext_{G_1}^1(\Delta(i), \Delta(p(a-b)+j)) 
                               &\mbox{if $a\ge b$}\\
                 \end{array} \right.}\\
            &\cong \left\{\begin{array}{ll}
                    \Delta(a-b-2)\frob &\mbox{if $a-b \ge 2$ and $i=j$}\\  
                    \nabla(b-a+1)\frob &\mbox{if $a-b \le 1$ and
                                                          $i=\bj$}\\  
                    0              &\mbox{otherwise}  
                 \end{array} \right.
\end{align*}
\end{lem}
\begin{proof}
We use the techniques of \cite{erd1}, where the lemma is proved for 
$a\ge b$.
Thus we need only prove the lemma for $b \ge a$.
We have 
\begin{align*}
\Ext^1_{G_1}(\Delta(pb+j), \Delta(pa+i)) 
&\cong
\Ext^1_{G_1}(\Delta(p(b-1)+\bj), \Delta(p(a-1)+\bi)) \\
&\cong
\Ext^1_{G_1}(\Delta(p(b-1)+j), \Delta(p(a-1)+i)) 
\end{align*}
for $a \ge1$ and $b\ge 1$, using the \ses \eqref{xanth2} twice.
The last isomorphism follows using the translation principle.
We thus get for $b \ge a$ that
$$
\Ext^1_{G_1}(\Delta(pb+j), \Delta(pa+i)) 
\cong
\Ext^1_{G_1}(\Delta(p(b-a)+j), \Delta(i))
$$
Now apply $\Hom_{G_1}(-,\Delta(i))$ to the \ses \eqref{xanth2} for
$\Delta(p(b-a)+j)$. We get
\begin{multline*}
0 \to
\Hom_{G_1}(\Delta(p(b-a)+j),\Delta(i))
\to
\Hom_{G_1}(\Delta(b-a)\frob \otimes Q(j),\Delta(i))\\
\to
\Hom_{G_1}(\Delta(p(b-a+1)+\bj),\Delta(i))
\to
\Ext^1_{G_1}(\Delta(p(b-a)+j),\Delta(i))
\to 0
\end{multline*}
The last zero follows as $Q(j)$ is projective. 
Since $\Delta(p(b-a)+j)$ has $G_1$-head $\Delta(b-a)\frob \otimes
L(j)$ we get that the first $\Hom$-group is $\nabla(b-a)\frob$ if
$j=i$ and $0$ otherwise. The second $\Hom$-group is also
$\nabla(b-a)\frob$ if $j=i$ and zero otherwise.
So
$$
\Hom_{G_1}(\Delta(p(b-a+1)+\bj),\Delta(i))
\cong
\Ext^1_{G_1}(\Delta(p(b-a)+j),\Delta(i)).
$$
This $\Hom$ group is $\nabla(b-a+1)\frob$ if $i=\bj$ and $0$
otherwise.
Hence the result.
\end{proof}

\begin{lem}\label{lem:nextone}
\begin{align*}
\Ext^1_{G_1}(\Delta(pb+i), \nabla(pa+j)) 
    &\cong \Ext_{G_1}^1(\Delta(p(a+b)+i), \Delta(j)) \\
    &\cong \left\{\begin{array}{ll}
            \nabla(a+b+1)\frob &\mbox{if $i=\bj$}\\  
            0              &\mbox{otherwise}  
     \end{array} \right.
\end{align*}
\end{lem}
\begin{proof}
A proof of this may be found in \cite[proposition 1.1]{devext1}
and is similar to that of the previous lemma.
\end{proof}

\begin{lem}\label{lem:Dhoms}
If $p\ge3$ then
$$
\Hom_{G_1}(\Delta(pb+\bi), \Delta(pa+i)) \cong
                    \Delta(a-1)\frob \otimes \nabla(b)\frob 
$$
(where $\Delta(-1)$ is interpreted as the zero module)
and
$$
\Hom_{G_1}(\Delta(pb+i), \Delta(pa+i)) 
          \cong \left\{\begin{array}{ll}
                    \nabla(b-a)\frob &\mbox{if $b\ge a$}\\
                    0              &\mbox{otherwise}  
                 \end{array} \right.
$$
If $p=2$ then
$$
\Hom_{G_1}(\Delta(2b), \Delta(2a)) 
          \cong \left\{\begin{array}{ll}
                    \Delta(a-1)\frob \otimes \nabla(b)\frob 
                    \oplus \nabla(b-a)\frob &\mbox{if $b\ge a$}\\
                    \Delta(a-1)\frob \otimes \nabla(b)\frob 
                                  &\mbox{if $b < a$}  
                 \end{array} \right.
$$
(where again $\Delta(-1)$ is interpreted as the zero module)
\end{lem}
\begin{proof}
We use the fact that the $G_1$-socle of $\Delta(pa+i)$ is
$\Delta(a-1)\frob \otimes L(\bi)$ if $a\ge1$
and the $G_1$-head of $\Delta(pa+i)$ is $\Delta(a) \frob \otimes L(i)$. 

We first assume that $p\ge3$ and so 
we may assume that $a\ge1$ as $\Delta(i)$ has the wrong $G_1$ type to
give a homomorphism.

Apply $\Hom_{G_1}(-,\Delta(pa+i))$ to the \ses \eqref{xanth1} for
$\Delta(pb+\bi)$.
\begin{multline*}
0 \to
\Hom_{G_1}(\Delta(b)\frob \otimes L(\bi),\Delta(pa+i))
\to
\Hom_{G_1}(\Delta(pb+\bi),\Delta(pa+i))\\
\to
\Hom_{G_1}(\Delta(b-1)\otimes L(i),\Delta(pa+i))
\end{multline*}
The first $\Hom$ group is $\nabla(b)\frob \otimes \Delta(a-1)\frob$.
The third $\Hom$ group is $0$. 
Hence the first result.

To prove the second we consider $\Hom_{G_1}(\Delta(pb+i),-)$ applied
to the \ses \eqref{xanth2} for $\Delta(a)\frob \otimes Q(i)$.
\begin{multline*}
0 \to
\Hom_{G_1}(\Delta(pb+i),\Delta(p(a+1)+\bi))
\to
\Hom_{G_1}(\Delta(pb+i),\Delta(a)\frob \otimes Q(i))\\
\to
\Hom_{G_1}(\Delta(pb+i),\Delta(pa+i))
\to
\Ext^1_{G_1}(\Delta(pb+i),\Delta(p(a+1)+\bi))
\to 0
\end{multline*}
The last zero follows as $Q(i)$ is injective. We use the first result
and lemma \ref{lem:Mfrob} to deduce that the first two $\Hom$ groups are
isomorphic to $\nabla(b)\frob \otimes \Delta(a)\frob$. 
Thus
$$
\Hom_{G_1}(\Delta(pb+i),\Delta(pa+i))
\cong
\Ext^1_{G_1}(\Delta(pb+i),\Delta(p(a+1)+\bi))
$$
This latter $\Ext$ group is $\nabla(b-a)\frob$ using lemma
\ref{lem:dextone}.

If $p=2$ then we prove the lemma by induction on $a$. 
It is clear for $a=0$ 
as we know that the $G_1$-head of $\Delta(2b)$ is $\Delta(b)\frob$.


We consider $\Hom_{G_1}(\Delta(2b),-)$ applied
to the \ses \eqref{xanth2} for $\Delta(a)\frob \otimes Q(0)$.
\begin{multline*}
0 \to
\Hom_{G_1}(\Delta(2b),\Delta(2(a+1)))
\to
\Hom_{G_1}(\Delta(2b),\Delta(a)\frob \otimes Q(0))\\
\to
\Hom_{G_1}(\Delta(2b),\Delta(2a))
\to
\Ext^1_{G_1}(\Delta(2b),\Delta(2(a+1)))
\to 0
\end{multline*}
The last zero follows as $Q(0)$ is injective. By induction and 
lemmas \ref{lem:Mfrob} and \ref{lem:dextone} we have that the 
first $\Hom$ group has the required character.
If $b \ge a$ then we have
\begin{multline*}
0 \to
\Hom_{G_1}(\Delta(2b),\Delta(2(a+1)))
\to
\nabla(b)\frob \otimes \Delta(a)\frob \oplus \nabla(b-1)\frob
\otimes \Delta(a)\frob\\
\stackrel{\phi}\to
\nabla(b)\frob \otimes \Delta(a-1)\frob \oplus 
\nabla(b-a)\frob
\to
\nabla(b-a)\frob
\to 0
\end{multline*}
Now if $b \ge a$ then
$ \nabla(b)\frob \otimes \Delta(a-1)\frob$ has a filtration by modules 
(in order starting at the top) $\nabla(a+b-1)\frob$,
$\nabla(a+b-3)\frob$, $\ldots$, $\nabla(b-a+1)\frob$ and
$\nabla(b-1)\frob \otimes \Delta(a)\frob$ has a filtration by modules 
(in order starting at the top) $\nabla(a+b-1)\frob$,
$\nabla(a+b-3)\frob$, $\ldots$, $\nabla(b-a-1)\frob$, using lemma
\ref{lem:good}. 
Also 
$ \nabla(b)\frob \otimes \Delta(a)\frob$ has a filtration by modules 
(in order starting at the top) $\nabla(a+b)\frob$,
$\nabla(a+b-2)\frob$, $\ldots$, $\nabla(b-a)\frob$ by lemma
\ref{lem:good}.
Note that this doesn't intersect with any factor in $\nabla(b)\frob
\otimes \Delta(a-1)\frob$.
Consequently the kernel of the map $\phi$ is
$\nabla(b)\frob \otimes \Delta(a)\frob \oplus \nabla(b-a-1)\frob$,
which is the required $\Hom$ group. 
 
We may do a similar argument if $b < a$ only now $\nabla(b)\frob
\otimes \Delta(a)\frob$ has a filtration by twisted Weyl modules.
%
\end{proof}

\begin{lem}\label{lem:nhoms}
$$
\Hom_{G_1}(\Delta(pb+i), \nabla(pa+j)) 
     \cong \left\{\begin{array}{ll}
         \nabla(b)\frob \otimes \nabla(a)\frob &\mbox{if $i=j$}\\
         0 &\mbox{otherwise}
      \end{array} \right.
$$
\end{lem}
\begin{proof}
A proof of this may be found in \cite[proposition 1.2]{devext1}
and is similar to that of the previous lemma.
%
%
\end{proof}

We finish this section by noting some vanishing results for the
$G$-cohomology.
\begin{propn}\label{propn:vanish}
Let $M$, $N \in \Mod(G)$ and set $pm_1+m_0$ to be the highest weight
of $M$ and $pn_1 +n_0$ to be the highest weight of $N$ with $m_0$,
$n_0\in X_1$ and $m$, $n \in X^+$, then
$$\Ext^i_{G}(M, N) = 0 \quad \mbox{if $i > m_1 +n_1$}$$
$$\Ext^i_G(M,\nabla(pn_1+n_0)) = 0\quad  \mbox{if $i > m_1$}$$
$$\Ext^i_G(\Delta(pm_1+m_0),N) = 0\quad  \mbox{if $i > n_1$}$$
\end{propn}
\begin{proof}
This follows using the results of \cite{parker1}, as $m_1$ is an upper
bound for the Weyl filtration dimension of $M$ and $n_1$ is an upper
bound for the good filtration dimension of $N$.
\end{proof}

\section{Applying to modules for $\SL_2(k)$ I}\label{sect:I}
We now apply the theory of section \ref{sect:ss} and look more closely
at the spectral sequences for $\Ext_G^q(V, \Delta(pa+i))$.

We first explicitly 
construct the $G_1$ injective resolution of $\Delta(pa+i)$.
This is easy using the sequence \eqref{xanth2}, and its dual, 
as the middle term is projective as a $G_1$-module.
We have 
$$ 
0\to 
\Delta(pa+i) 
\xrightarrow{\delta_{-1}}
I_0 \xrightarrow{\, \delta_0\, }
I_1 \xrightarrow{\, \delta_1\,}
\cdots  \xrightarrow{\delta_{m-1}}
I_m  \xrightarrow{\delta_{m}}
I_{m+1} \xrightarrow{\delta_{m+1}} \cdots  
$$
with   
$$
I_m \cong \left\{\begin{array}{ll}
                    \Delta(a-m-1)\frob \otimes Q(i) &\mbox{if $m$ odd
                                  and $m \le a-1$}\\
                    \Delta(a-m-1)\frob \otimes Q(\bi) &\mbox{if $m$
                                  even and $m \le a-1$}\\
                    \nabla(m-a)\frob \otimes Q(\bi) &\mbox{if $m$ odd
                                  and $m \ge a$}\\
                    \nabla(m-a)\frob \otimes Q(i) &\mbox{if $m$
                                  even and $m \ge a$}
                 \end{array} \right.
$$
We define $M_m = \ker\delta_{m}$. We have
$$
M_m \cong \left\{\begin{array}{ll}
                    \Delta(p(a-m)+\bi)&\mbox{if $m$ odd
                                  and $m \le a-1$}\\
                    \Delta(p(a-m)+i) &\mbox{if $m$
                                  even and $m \le a-1$}\\
                    \nabla(p(m-a)+\bi)   &\mbox{if $m$ odd
                                  and $m \ge a$}\\
                    \nabla(p(m-a)+i)   &\mbox{if $m$
                                  even and $m \ge a$}
                 \end{array} \right.
$$

We have for $N$ a $G$-module and $p\ge3$ that
$$
\Hom_{G_1}(N\frob \otimes L(i), I_m) 
          \cong \left\{\begin{array}{ll}
                    \Delta(a-m-1)\frob \otimes (N\frob)^* &\mbox{if $m$ odd
                                  and $m \le a-1$}\\
                    \nabla(m-a)\frob \otimes (N\frob)^* &\mbox{if $m$
                                  even and $m \ge a$}\\
                    0              &\mbox{otherwise}  
                 \end{array} \right.
$$
and
$$
\Hom_{G_1}(N\frob \otimes L(\bi), I_m) 
          \cong \left\{\begin{array}{ll}
                    \Delta(a-m-1)\frob \otimes (N\frob)^* &\mbox{if $m$
                                  even and $m \le a-1$}\\
                    \nabla(m-a)\frob \otimes (N\frob)^* &\mbox{if $m$ odd
                                  and $m \ge a$}\\
                    0              &\mbox{otherwise.}  
                 \end{array} \right.
$$
If $p=2$ then
$$
\Hom_{G_1}(N\frob, I_m) 
          \cong \left\{\begin{array}{ll}
                    \Delta(a-m-1)\frob \otimes (N\frob)^* &\mbox{if $m \le a-1$}\\
                    \nabla(m-a)\frob \otimes (N\frob)^* &\mbox{if $m
                      \ge a$.}\\
                 \end{array} \right.
$$
If $p\ge3$ then we also have 
$$
\Hom_{G_1}(N\frob \otimes L(i), M_m) 
          \cong \left\{\begin{array}{ll}
                    \Delta(a-m-1)\frob \otimes (N\frob)^* &\mbox{if $m$ odd
                                  and $m \le a-1$}\\
                    \nabla(m-a)\frob \otimes (N\frob)^* &\mbox{if $m$
                                  even and $m \ge a$}\\
                    0              &\mbox{otherwise}  
                 \end{array} \right.
$$
and
$$
\Hom_{G_1}(N\frob \otimes L(\bi), M_m) 
          \cong \left\{\begin{array}{ll}
                    \Delta(a-m-1)\frob \otimes (N\frob)^* &\mbox{if $m$
                      even
                                  and $m \le a-1$}\\
                    \nabla(m-a)\frob \otimes (N\frob)^* &\mbox{if $m$
                                  odd and $m \ge a$}\\
                    0              &\mbox{otherwise.}  
                 \end{array} \right.
$$

If $p=2$ then
$$
\Hom_{G_1}(N\frob, M_m) 
          \cong \left\{\begin{array}{ll}
                    \Delta(a-m-1)\frob \otimes (N\frob)^* &\mbox{if $m \le a-1$}\\
                    \nabla(m-a)\frob \otimes (N\frob)^* &\mbox{if $m
                      \ge a$}\\
                 \end{array} \right.
$$
using lemmas \ref{lem:Dhoms} and \ref{lem:nhoms}. 

We have exact sequences
\begin{multline*}
0 \to
\Hom_{G_1}(N\frob\otimes L(j), M_m) \to
\Hom_{G_1}(N\frob\otimes L(j), I_m) \\
 \xrightarrow{\delta_{m}^*}
\Hom_{G_1}(N\frob\otimes L(j), M_{m+1}) \to
\Ext_{G_1}^1(N\frob\otimes L(j), M_m) \to
0.
\end{multline*}

Now since we have
$$
\Hom_{G_1}(N\frob \otimes L(j), I_m) 
\cong
\Hom_{G_1}(N\frob \otimes L(j), M_m) 
$$
for all $m$, $j$ and $p$, 
we must have
\begin{align*}
\Ext^m_{G_1}(N\frob \otimes L(j), \Delta(pa+i)) 
&\cong \Ext^1_{G_1}(N\frob \otimes L(j), M_{m-1}) \\
&\cong \Hom_{G_1}(N\frob \otimes L(j), M_m) \\
&\cong\Hom_{G_1}(N\frob \otimes L(j), I_m). 
\end{align*}

Thus
all the induced maps $\delta^*_m$ in the
following complex are zero.
\begin{multline*}
\Hom_{G_1}(N\frob\otimes L(j), I_0) \xrightarrow{\delta_0^*}
\Hom_{G_1}(N\frob\otimes L(j), I_1) \xrightarrow{\delta_1^*}
\Hom_{G_1}(N\frob\otimes L(j), I_2) 
\\  \xrightarrow{\delta_2^*}
\cdots
\xrightarrow{\delta_{m-2}^*}
\Hom_{G_1}(N\frob\otimes L(j), I_{m-1}) \xrightarrow{\delta_{m-1}^*}
\Hom_{G_1}(N\frob\otimes L(j), I_{m}) 
\\ \xrightarrow{\delta_{m}^*}
\Hom_{G_1}(N\frob\otimes L(j), I_{m+1}) \xrightarrow{\delta_{m+1}^*}
\end{multline*}

Since all the differentials in the $\Hom_{G_1}$ complex are zero, the
induced differentials $d_0$ on the $E_0$ page are also zero.
Thus  all the $k_2$'s associated to the 2nd derived couple
are zero using lemma \ref{lem:dzero}. 
Hence the $E_2$ page is the same as the $E_\infty$ page,
using corollary \ref{cor:E2}.

Now the $E_2$ page for the spectral sequence associated to
$\Ext^q_G(N\frob \otimes L(i), \Delta(pa+i))$ and $p\ge3$
is
\begin{align*}
E_2^{m,n}
&= H^m(G/G_1, \Ext^n_{G_1}(N\frob \otimes L(i), \Delta(pa+i)))\\
&\cong \left\{\begin{array}{ll}
                    \Ext^m_G(N, \Delta(a-n-1))&\mbox{if $n$ odd
                                  and $n \le a-1$}\\
                    \Ext^m_G(N, \nabla(n-a)) &\mbox{if $n$
                                  even and $n \ge a$}\\
                    0              &\mbox{otherwise}  
                 \end{array} \right.
\end{align*}
and for $\Ext^q_G(\nabla(pa+i), N\frob \otimes L(\bi))$ and $p\ge3$ we have
\begin{align*}
E_2^{m,n}
&= H^m(G/G_1, \Ext^n_{G_1}(\nabla(pa+i), N\frob \otimes L(\bi)))\\
&\cong \left\{\begin{array}{ll}
                    \Ext_G^m(N, \Delta(a-n-1)) &\mbox{if $n$
                                  even and $n \le a-1$}\\
                    \Ext^m_G(N, \nabla(n-a))&\mbox{if $n$
                                  odd and $n \ge a$}\\
                    0              &\mbox{otherwise}  
                 \end{array} \right.
\end{align*}

If $p=2$ then 
the $E_2$ page for the spectral sequence associated to
$\Ext^q_G(N\frob, \Delta(2a))$ is
\begin{align*}
E_2^{m,n}
&= H^m(G/G_1, \Ext^n_{G_1}(N\frob, \Delta(2a)))\\
&\cong \left\{\begin{array}{ll}
                    \Ext^m_G(N, \Delta(a-n-1))&\mbox{if $n \le a-1$}\\
                    \Ext^m_G(N, \nabla(n-a)) &\mbox{if $n \ge a$.}\\
                 \end{array} \right.
\end{align*}

We thus have the following theorem.
\begin{thm} \label{thm:sectI}
Let $a$, $q \in \N$, $0 \le i \le p-2$,  and $N \in \Mod(G)$.
If $p\ge 3$ we have
\begin{align*}
\Ext^q_{G}(N\frob\otimes L(i), &\ \Delta(pa+i))\\
&\cong \bigoplus_{\substack{n\mbox{\rm{\scriptsize{ odd}}}\\ \\{0 \le
      n \le \min\{q, a-1\}}}}
      \Ext^{q-n}_G(N, \Delta(a-n-1))
 \oplus \bigoplus_{\substack{n\mbox{\rm{\scriptsize{ even}}} 
                       \\ \\ a \le n \le q}}
      \Ext^{q-n}_G(N,\nabla(n-a)) 
\end{align*}
and
\begin{align*}
\Ext^q_{G}(N\frob\otimes L(\bi), &\ \Delta(pa+i))\\
&\cong \bigoplus_{\substack{n \mbox{\rm{\scriptsize{ even}}}\\ \\{0 \le
      n \le \min\{q, a-1\}}}}
      \Ext^{q-n}_G(N, \Delta(a-n-1))
\oplus \bigoplus_{\substack{n\mbox{\rm{\scriptsize{ odd}}} 
                     \\ \\ a \le n \le q}}
      \Ext^{q-n}_G(N,\nabla(n-a)) 
\end{align*}
If $p=2$ then
$$
 \Ext^q_{G}(N\frob, \Delta(2a))) 
\cong
  \bigoplus_{n=0}^{n=\min\{q,a-1\}}
       \Ext^{q-n}_G(N, \Delta(a-n-1))
  \oplus \bigoplus_{n=a}^{n=q}
       \Ext^{q-n}_G(N, \nabla(n-a)) 
$$
\end{thm}
\begin{cor}
Let $a \ge 1$, $0 \le i \le p-2$, $q \in \N$ and $N \in \Mod(G)$.
If $p\ge3$ we have
$$
\Ext^q_{G}(N\frob\otimes L(i), \Delta(pa+i)) 
\cong \Ext_{G}^{q-1}(N\frob\otimes L(i), 
\Delta(p(a-1)+\bi))
$$
$$
\Ext_{G}^{q}(N\frob\otimes L(i), \Delta(pa +\bi))
\cong \Ext_{G}^{q-1}(N\frob\otimes L(i), 
\Delta(p(a-1)+i))
\oplus \Ext_{G}^{q}(N, \Delta(a-1))
$$
where $\Ext^{-1}$ is interpreted as the zero module.

If $p=2$ then
$$
 \Ext^q_{G}(N\frob, \Delta(2a))) 
\cong
       \Ext^{q-1}_G(N\frob, \Delta(2a-2))
  \oplus \Ext^{q}_G(N, \Delta(a-1))
$$
\end{cor}
Hence if $N=L(b)$ then this would completely determine
$\Ext^q_G(L(pb+j), \Delta(pa+i))$ by induction, once we knew what
$\Ext^q_G(L(b), \nabla(n-a))$ was.
Similarly for $N=\nabla(b)$ or $\Delta(b)$. (This is as
$\Ext^q_G(N,\Delta(a))$ will be zero if $q$ is larger than the highest
weight of $N$ divided by $p$ plus $\lfloor\frac{a}{p}\rfloor$ using lemma
  \ref{propn:vanish}.)

If $N=k$ then the $\Ext^m_G(k, \nabla(n-a))$ vanish unless $n=a$ and
$m=0$, so this is enough to work out $\Ext^q_G(\Delta(j), \Delta(pa+i))$. 
This particular $\Ext$ group is also a special case of the results in
section \ref{sect:II}.

We may now use theorem \ref{thm:sectI}
to deduce the following theorem where we set $a=0$, ``$N$'' to be
$N\otimes M^*$ and ``$M$'' to be $k$.
\begin{thm}\label{thm:sectI2}
Let $N$ and $M$ be in $\Mod(G)$ and $0 \le i \le p-2$.
Then if $p\ge 3$
$$
 \Ext^q_{G}(N\frob\otimes L(i), M\frob \otimes L(i))) 
\cong \bigoplus_{\substack{n \mbox{\rm{\scriptsize{ even}}}\\ \\{0 \le
      n \le q}}}
      \Ext^{q-n}_G(N, \nabla(n)\otimes M)
$$
and
$$
\Ext^q_{G}(N\frob\otimes L(i), M\frob \otimes L(\bi))) 
\cong \bigoplus_{\substack{n \mbox{\rm{\scriptsize{ odd}}}\\ \\{0 \le
      n \le q}}}
      \Ext^{q-n}_G(N, \nabla(n)\otimes M)
$$
If $p=2$ then
$$
 \Ext^q_{G}(N\frob, M\frob) 
\cong
\bigoplus_{n=0}^{q}
     \Ext^{q-n}_G(N, \nabla(n)\otimes M).
$$
\end{thm}

So this would completely determine $\Ext^q_G(L(pa+i), L(pb+j)$ if we
knew what $\Ext^{q}_G (L(a), \nabla(n)\break \otimes L(b))$ was.
This in principle we could calculate using the following:
$$\Ext^{q}_G (L(a), \nabla(n)\otimes L(b))
\cong
\Ext^{q}_G (L(a)\otimes L(b), \nabla(n))$$
since $L(b)^* \cong L(b)$. The structure of the tensor product 
$L(a)\otimes L(b)$ has
been determined in \cite{dothen}. We get a direct sum of
modules of the form $M\frob \otimes T(l)$ where $l$ is either
$p$-restricted, in which case we could use theorem \ref{thm:DL} to
calculate the $\Ext$ group or $p-1 \le l \le 2p-2$ in which case we
would use the dual version of lemma \ref{lem:MT} to calculate the
$\Ext$ group.

We can get a complete formula if $N=\Delta(a)$ and $M=\nabla(b)$.
Here $\Ext^{m}_G (\Delta(a), \nabla(n) \otimes \nabla(b))$ is
zero if $m \ge 1$,
and $\Hom_G(\Delta(a), \nabla(n)\otimes \nabla(b))$ is either $0$ or
$k$ and it is $k$ if and only if
$\nabla(a)$ is a section of $\nabla(n)\otimes \nabla(b)$. 
This happens if $a$ is of the same parity as $n+b$ 
and $n+b \ge a \ge \max\{n-b, b-n\}$ 

Thus we get 
$$
\Ext^q_{G}(\Delta(a)\frob\otimes L(i), \nabla(b)\frob \otimes L(j))) 
\cong
\left\{\begin{array}{ll}
             k &\mbox{if $q$ even, $a+b$ even,  $i=j$}\\ 
                  &\mbox{\hspace{20pt}and $q+b \ge a \ge \max\{q-b, b-q\}$ }\\
             k &\mbox{if $q$ odd, $a+b$ odd, $i=\bj$}\\ 
                  &\mbox{\hspace{20pt}and $q+b \ge a \ge \max\{q-b, b-q\}$ }\\
             0 &\mbox{otherwise}
\end{array} \right.
$$

We consider the special case where $pa+i$ and $pb+j$ are in 
the Jantzen region and
so both $a$ and $b$ are at most $p-1$. We have $L(pa+i) \cong
\Delta(a)\frob\otimes L(i)$ and
$L(pb+j) \cong \nabla(b)\frob\otimes L(j)$ this then gives the 
$\Ext$ between simples and hence the coefficients of the
dual Khazdhan-Lusztig polynomials.

\section{Applying to modules for $\SL_2(k)$ part II}\label{sect:II}

We now do a similar procedure ---
we still use the $G_1$-injective resolution for $\Delta(pa+i)$
constructed in the previous section but now we take $V=\Delta(pb+j)$.

We first assume that $p\ge 3$. We have
$$
\Hom_{G_1}(\Delta(pb+i), I_m) 
          \cong \left\{\begin{array}{ll}
                    \Delta(a-m-1)\frob \otimes \nabla(b)\frob &\mbox{if $m$ odd
                                  and $m \le a-1$}\\
                    \Delta(a-m-1)\frob \otimes \nabla(b-1)\frob
                                  &\mbox{if $m$ even
                                  and $m \le a-1$}\\
                    \nabla(m-a)\frob \otimes \nabla(b)\frob &\mbox{if $m$
                                  even and $m \ge a$}\\
                    \nabla(m-a)\frob \otimes \nabla(b-1)\frob &\mbox{if $m$
                                  odd and $m \ge a$}\\
                 \end{array} \right.
$$
and
$$
\Hom_{G_1}(\Delta(pb+\bi), I_m) 
          \cong \left\{\begin{array}{ll}
                    \Delta(a-m-1)\frob \otimes \nabla(b-1)\frob &\mbox{if $m$ odd
                                  and $m \le a-1$}\\
                    \Delta(a-m-1)\frob \otimes \nabla(b)\frob
                                  &\mbox{if $m$ even
                                  and $m \le a-1$}\\
                    \nabla(m-a)\frob \otimes \nabla(b-1)\frob &\mbox{if $m$
                                  even and $m \ge a$}\\
                    \nabla(m-a)\frob \otimes \nabla(b)\frob &\mbox{if $m$
                                  odd and $m \ge a$}\\
                 \end{array} \right.
$$

So we have using lemmas \ref{lem:Dhoms} and \ref{lem:nhoms} that
$$
\Hom_{G_1}(\Delta(pb+i), M_m) 
          \cong \left\{\begin{array}{ll}
                    \Delta(a-m-1)\frob \otimes \nabla(b)\frob &\mbox{if $m$ odd
                                  and $m \le a-1$}\\
                    \nabla(m-a+b)\frob &\mbox{if $m$ even} \\
                                &\mbox{\hspace{20pt}and $a-b \le m\le a-1$}\\ 
                    \nabla(m-a)\frob \otimes \nabla(b)\frob &\mbox{if
                                   $m$ even and $m \ge a$}\\
                    0              &\mbox{otherwise}  
                 \end{array} \right.
$$
and
$$
\Hom_{G_1}(\Delta(pb+\bi), M_m) 
          \cong \left\{\begin{array}{ll}
                    \Delta(a-m-1)\frob \otimes \nabla(b)\frob &\mbox{if
                                 $m$ even and $m \le a-1$}\\
                    \nabla(m-a+b)\frob &\mbox{if $m$ odd}\\
                                &\mbox{\hspace{20pt}and $a-b \le m\le a-1$}\\ 
                    \nabla(m-a)\frob \otimes \nabla(b)\frob &\mbox{if
                                   $m$ odd and $m \ge a$}\\
                    0              &\mbox{otherwise}  
                 \end{array} \right.
$$

So 
$$
\Hom_{G_1}(\Delta(pb+i), I_m) 
\cong
\Hom_{G_1}(\Delta(pb+i), M_m) $$
if $m$ odd and $m \le a-1$ or if $m$ even and $m \ge a$
and 
$$
\Hom_{G_1}(\Delta(pb+\bi), I_m) 
\cong
\Hom_{G_1}(\Delta(pb+\bi), M_m) $$
if $m$ even and $m \le a-1$ or if $m$ odd and $m \ge a$

Thus the induced differential $\delta_{m}^*$ on the $\Hom_{G_1}$
complex is zero in the above cases.

We have exact sequences
\begin{multline*}
0 \to
\Hom_{G_1}(\Delta(pb+j), M_m) \to
\Hom_{G_1}(\Delta(pb+j), I_m) \\\xrightarrow{\delta_{m}^*}
\Hom_{G_1}(\Delta(pb+j), M_{m+1}) \to
\Ext_{G_1}^1(\Delta(pb+j), M_m) \to
0
\end{multline*}

We thus get for $m\ge 1$ using lemmas \ref{lem:dextone} and
\ref{lem:nextone}
\begin{align*}
\Ext^m_{G_1}(\Delta(pb+i), \Delta(pa+i)) 
          &\cong \Ext_{G_1}^1(\Delta(pb+i), M_{m-1}) \\
&\cong \left\{\begin{array}{ll}
                    \Delta(a-m-b-1)\frob &\mbox{if $m$ odd}\\  
                                &\mbox{\hspace{35pt}and $m \le a-b-1$}\\
                    \nabla(m-a+b)\frob &\mbox{if $m$ even}\\
                                &\mbox{\hspace{35pt}and $m\ge a-b$}\\
                    0              &\mbox{otherwise}  
                 \end{array} \right.
\end{align*}
and
\begin{align*}
\Ext^m_{G_1}(\Delta(pb+\bi), \Delta(pa+i)) 
          &\cong \Ext_{G_1}^1(\Delta(pb+\bi), M_{m-1}) \\
&\cong \left\{\begin{array}{ll}
                    \Delta(a-m-b-1)\frob &\mbox{if $m$ even}\\  
                                &\mbox{\hspace{35pt}and $m \le a-b-1$}\\
                    \nabla(m-a+b)\frob &\mbox{if $m$ odd}\\
                                &\mbox{\hspace{35pt}and $m\ge a-b$}\\
                    0              &\mbox{otherwise.}  
                 \end{array} \right.
\end{align*}

We now consider the $E_1$ page corresponding to $\Ext^q_G(\Delta(pb+i),
\Delta(pa+j))$.
As in the previous section, we will show that the $k_2$'s are all
zero.

We consider the case for $\Ext^q_G(\Delta(pb+i), \Delta(pa+i))$, the
other case is similar.
We have that the differential $d_0: E_0^{m,n}\to E_0^{m,n+1}$ is zero
for $n$ odd and $n \le a-1$ or for $n$ even and $n \ge a$.
We also know that $E_1^{m,n+1} = \Hom_{G/{G_1}}(k,
\Ext^{n+1}_{G_1}(\Delta(pb+i), \Delta(pa+i))\otimes J_m\frob)$
 is zero for $n$ odd and 
$n\le a-b-2$ or for $n$ even and $n \ge a-b-1$.
But since $E_1^{m,n+1}$ is the homology at $E_0^{m,n+1}$ with respect
to $d_0$ we must have that the kernel of $d_0:E_0^{m,n+1}\to
E_0^{m,n+2}$
is the image of $d_0:E_0^{m,n}\to E_0^{m,n+1}$ which is zero for $n$
odd and $n\le a-b-2< a-1$.
Thus the map $d_0:E_0^{m,n+1}\to E_0^{m,n+2}$ is an embedding, for
$n$ odd and $n\le a-b-2$. Note this is independent of $m$.

Thus for $n$ odd and $n \le a-b-2$ we have that 
the map $d_0: E_0^{m,n}\to E_0^{m,n+1}$ is zero 
and that the map $d_0: E_0^{m+2,n-1} \to E_0^{m+2,n}$ is an embedding. 
Now by lemmas \ref{lem:dzero} and \ref{lem:dinj} the maps 
$k_2:E_2^{m-1,n+1} \to D_2^{m,n+1}$ and
$k_2:E_2^{m,n} \to D_2^{m+1,n}$ 
are zero.

Now $k_2^{m,n}$ is of course zero if $E_2^{m,n}$ is.
We have 
\begin{align*}
E_2^{m,n} &= H^m(G/G_1, \Ext^n_{G_1}(\Delta(pb+i), \Delta(pa+i)) \\
&= \left\{\begin{array}{ll}
                    \Ext^m_G(k,\Delta(a-n-b-1)) 
                                &\mbox{if $n$ odd and $n \le a-b-1$}\\
                    \Ext^m_G(k,\nabla(n-a+b)) 
                                &\mbox{if $n$ even and $n\ge a-b$}\\
                    0              &\mbox{otherwise.}  
                 \end{array} \right.\\
&= \left\{\begin{array}{ll}
                    \Ext^m_G(k, \nabla(a-n-b-1)) 
                                &\mbox{if $n$ odd and $n \le a-b-1$}\\
                    k &\mbox{if $n$ even, $n=a-b$ and $m=0$}\\
                    0              &\mbox{otherwise.}  
                 \end{array} \right.
\end{align*}

Thus the only $k_2$ left to show is zero is
$k_2^{0,a-b}: E_2^{0,a-b} \to D_2^{1,a-b}$, but now $d_0: E_0^{1,a-b-1}
\to E_0^{1,a-b}$ is zero and so lemma \ref{lem:dzero}
shows that $k_2^{0,a-b}$ is zero.

Hence all the $k_2$'s in the second derived couple are zero. This in
particular implies that all subsequent differentials $d_f$ are zero
for $f \ge 2$.  Thus the $E_2$ page is the same as the $E_\infty$
page.

We may argue in exactly the same way to show that the $k_2$'s for the
second derived couple for
$\Ext^q_G(\Delta(pb+\bi), \Delta(pa+i))$ are also zero and that
the $E_2$ page is the same as the $E_\infty$ page.

We have 
\begin{align*}
E_2^{m,n} &= H^m(G/G_1, \Ext^n_{G_1}(\Delta(pb+\bi), \Delta(pa+i)) \\
&= \left\{\begin{array}{ll}
                    \Ext^m_G(k,\Delta(a-n-b-1)) 
                           &\mbox{if $n$ even and  $0 < n \le a-b-1$}\\
                    \Ext^m_G(\Delta(b),\Delta(a-1)) 
                                &\mbox{if $n=0 \le a-b-1$}\\
                    \Ext^m_G(k,\nabla(n-a+b)) 
                                &\mbox{if $n$ odd and $n\ge a-b$}\\
                    0              &\mbox{otherwise.}  
                 \end{array} \right.\\
&= \left\{\begin{array}{ll}
                    \Ext^m_G(k,\Delta(a-n-b-1)) 
                                &\mbox{if $n$ even and $0 < n \le a-b-1$}\\
                    \Ext^m_G(\Delta(b),\Delta(a-1)) 
                                &\mbox{if $n=0 \le a-b-1$}\\
                    k &\mbox{if $n$ even, $n=a-b$ and $m=0$}\\
                    0              &\mbox{otherwise.}  
                 \end{array} \right.
\end{align*}
We may thus deduce theorem \ref{thm:extD} for $p\ge3$.

We now assume that $p=2$. We have
$$
\Hom_{G_1}(\Delta(2b), I_m) 
    \cong \left\{\begin{array}{ll}
        \Delta(a-m-1)\frob \otimes \bigl(\nabla(b-1)\frob \oplus
        \nabla(b)\frob \bigr)
                &\mbox{if $m \le a-1$}\\
        \nabla(m-a)\frob \otimes \bigl(\nabla(b-1)\frob \oplus
        \nabla(b)\frob \bigr)
                &\mbox{if $m \ge a$.}\\
          \end{array} \right.
$$

We also have using lemmas \ref{lem:Dhoms} and \ref{lem:nhoms} that
$$
\Hom_{G_1}(\Delta(2b), M_m) 
    \cong \left\{\begin{array}{ll}
              \Delta(a-m-1)\frob \otimes \nabla(b)\frob 
                    &\mbox{if $m < a-b$}\\
              \Delta(a-m-1)\frob \otimes \nabla(b)\frob 
                    &\mbox{if $a-b \le m\le a-1$}\\ 
              \hspace{20pt} \oplus \nabla(m-a+b)\frob &\\
              \nabla(m-a)\frob \otimes \nabla(b)\frob 
                    &\mbox{if $m \ge a$.}\\
       \end{array} \right.
$$

We have exact sequences
\begin{multline*}
0 \to
\Hom_{G_1}(\Delta(2b), M_m) \xrightarrow{\phi_{m}}
\Hom_{G_1}(\Delta(2b), I_m) \\\xrightarrow{\delta_{m}^*}
\Hom_{G_1}(\Delta(2b), M_{m+1}) \to
\Ext_{G_1}^1(\Delta(2b), M_m) \to
0.
\end{multline*}

We assume that $a > b$ and that $m < a-b$
and we write $\Hom_{G_1}(\Delta(2b), I_m)$ as $A_m \oplus B_m$ where
$A_m \cong \Delta(a-m-1)\frob \otimes \nabla(b)\frob$ and
$B_m \cong \Delta(a-m-1)\frob \otimes \nabla(b-1)\frob$. 
Note that $\Hom_{G_1}(\Delta(2b), M_m) \cong A_m$.

We claim that the map $\delta_m^*$ is the natural projection onto 
$B_m$ followed
by the natural embedding into $\Hom_{G_1}(\Delta(2b), M_{m+1})$.

This is clear by considering the $G/G_1$-module structure of $A_m$ and
$B_m$. Since $b < a-m$ the module $A_m$ has filtration by
twisted Weyl modules, namely  (starting at the bottom) $\Delta(a+b-m-1)\frob$, 
$\Delta(a+b-m-3)\frob$, $\Delta(a+b-m-5)\frob$, $\ldots$,
$\Delta(a-b-m-1)\frob$.
The module $B_m$ also has a filtration by
twisted Weyl modules, namely (starting at the bottom) $\Delta(a+b-m-2)\frob$, 
$\Delta(a+b-m-4)\frob$, $\Delta(a+b-m-6)\frob$, $\ldots$,
$\Delta(a-b-m)\frob$.
Thus the parities of the highest weights of the twisted Weyl modules
are distinct. In particular they are in different blocks of $G/G_1$.
Hence the map $\phi_{m}: \Hom_{G_1}(\Delta(2b),
M_m) \to \Hom_{G_1}(\Delta(2b), I_m)$ is the natural injection into
$A_m$ and has image $A_m$ (as it is a $G/G_1$-module homomorphism).
So $A_m$ is the kernel of $\delta_m^*$ and
the map $\delta_m^*$ is as described.

We now claim that the induced map $d^{n,m}_0$ on the $E_0$ page is zero when
$a-b$ is even, $m$ is odd and $m < a-b$, 
or if
$a-b$ is odd, $m$ is even and $m < a-b$. 

Now $E^{n,m}_0 = \Hom_{G/G_1}(k, A_m\otimes J_n\frob 
\oplus B_m \otimes J_n\frob )$.
The map $d_0^{n,m}$ takes a homomorphism $\phi \in E_0^{n,m}$ to the
homomorphism
$(\delta_m^* \otimes \id) \circ \phi$ where $\id$ is
the identity map on $J_n\frob$.
So the image of $d^{n,m}_0$ is contained in $\Hom_{G/G_1}(B_m \otimes
J_n\frob)$ if $m < a-b$.

We have 
$$
\Hom_{G/G_1}(k, B_m \otimes J_n\frob)
\cong \Hom_{G/G_1}(B_m^*, J_n\frob) 
\cong 0
$$
if $a-b-m$ is odd and $m < a-b$ 
as then we know that the modules appearing in the filtration of
$B_m^*$ have the wrong highest weight to give a $G/G_1$-homomorphism to 
$J_n\frob$.
(We may assume that the components of $J_n\frob$ are all in the block
of $k$.)
Thus the map $d_0^{n,m}$ is zero 
if $a-b$ is even and $m < a-b$ is odd
or if 
$a-b$ is odd, $m$ is even and $m < a-b$, 

We get for $m\ge 1$ using lemmas \ref{lem:dextone} and
\ref{lem:nextone}
\begin{align*}
\Ext^m_{G_1}(\Delta(2b), \Delta(2a)) 
          &\cong \Ext_{G_1}^1(\Delta(2b), M_{m-1}) \\
&\cong \left\{\begin{array}{ll}
                    \Delta(a-m-b-1)\frob 
                                &\mbox{if $m \le a-b-1$}\\
                    \nabla(m-a+b)\frob 
                                &\mbox{if $m\ge a-b$.}\\
                 \end{array} \right.
\end{align*}

Note for $m +1< a-b$ that
\begin{align*}
E_1^{n,m+1} &= \Hom_{G/G_1}(k, \Delta(a-m-b-2)\frob \otimes J_n\frob)\\
 &\cong \Hom_{G/G_1}(\nabla(a-m-b-2)\frob, J_n\frob)
\end{align*}
and that this is zero if $a-m-b-2$ is odd.
So it is zero if $a-b$ is even and $m$ is odd or if
$a-b$ is odd and $m$ is even.

Now the same argument as for $p \ge 3$ allows us to deduce that
the map $d_0^{n,m+1}$ is an embedding, if $a-m-b-2$ is odd.
Thus the argument then proceeds as in the $p \ge 3$ case and
the $E_2$ page is the same as the $E_\infty$ page.

Suppose $a > b$ then
we have 
\begin{align*}
E&_2^{m,n} = H^m(G/G_1, \Ext^n_{G_1}(\Delta(2b), \Delta(2a)) \\
&= \left\{\begin{array}{ll}
        \Ext_G^m(\Delta(b), \Delta(a-1)) 
              &\mbox{if $n=0$}\\
        \Ext^m_G(k,\Delta(a-n-b-1)) 
              &\mbox{if $1 \le n \le a-b-1$}\\
        \Ext^m_G(k,\nabla(n-a+b)) 
              &\mbox{if $n\ge a-b \ge 1$}\\
      \end{array} \right.\\
&= \left\{\begin{array}{ll}
         \Ext_G^m(\Delta(b), \Delta(a-1)) 
              &\mbox{if $n=0$}\\
         \Ext^m_G(k, \Delta(a-n-b-1)) 
             &\mbox{if $a-b-n-1$ is even and $1 \le n \le a-b-1$}\\
         k &\mbox{if $n=a-b$ and $m=0$}\\
         0              &\mbox{otherwise}  
    \end{array} \right.
\end{align*}
since the weight $a-n-b-1$ is only in the same block as the weight $0$
if it is even.

We thus deduce the following theorem.
\begin{thm}\label{thm:extD}
Suppose $a-b$ is even then
\begin{multline*}
 \Ext^q_{G}(\Delta(pb+i), \Delta(pa+i)) \\
\cong
\left\{\begin{array}{ll}
                 \bigoplus_{n=0}^{n=(a-b-2)/2}
                    \Ext^{q-2n-1}_G(k,\Delta(a-2n-b-2)) 
                                &\mbox{if $q \le a-b-1$}\\
                    k &\mbox{if $q=a-b$ }\\
                    0              &\mbox{otherwise.}  
                 \end{array} \right.
\end{multline*}

Suppose $a-b$ is odd then
\begin{multline*}
 \Ext^q_{G}(\Delta(pb+\bi), \Delta(pa+i)) \\
\cong
\left\{\begin{array}{ll}
     \Ext^q_G(\Delta(b), \Delta(a-1)) &\mbox{if $q \le a-b-1$}\\
         \hspace{20pt}       \oplus \bigoplus_{n=0}^{n=(a-b-1)/2} 
                            \Ext^{q-2n-2}_G(k, \Delta(a-2n-b-3)) 
                                & \\
     k              &\mbox{if $q=a-b$}\\
     0              &\mbox{otherwise.}  
\end{array} \right.
\end{multline*}

\end{thm}
(Remark: this is true for $q=0,1,2$ by the results of \cite{coxerd,erd1})

\begin{cor}\label{cor:extD}
For $a \ge b$ and $a-b$ odd we have
\begin{align*}
\Ext_G^m(\Delta(pb+\bi),\  &\Delta(pa+i))\\
&\cong \Ext_G^{m-1}(\Delta(pb+\bi), \Delta(p(a-1) +\bi))
\oplus \Ext_G^m(\Delta(b), \Delta(a-1))
\end{align*}
For $a \ge b$ and $a-b$ even we have
\begin{align*}
\Ext_G^{m}(\Delta(pb+i), \Delta(pa +i))
&\cong \Ext_G^{m}(\Delta(i), \Delta(p(a-b)+i))\\
&\cong \Ext_G^{m-1}(\Delta(i), \Delta(p(a-b-1) +\bi))
\end{align*}
where the last equality follows if $m\ge1$.
\end{cor}
Thus $\Ext^q_G(\Delta(pa+i), \Delta(pb+j))$ may now be determined by
induction as $\Ext^q_G(\Delta(b), \Delta(a-1))$ will be zero if $q >
a-b-1$ using \cite{ryom}.

It is also now possible to completely determine
$\Ext^m_G(T(pb+j), \Delta(pa+i))$.
Since $T(pb+j) \cong T(b-1)\frob \otimes T(p+j)$
we have
\begin{multline*}
\Ext^m_G(T(pb+j), \Delta(pa+i))
\cong
\Ext_G^{m}(T(b-1)\frob \otimes T(p+j), \Delta(pa+i))\\
\cong  \Ext^m_G(T(b-1), \Delta(a-1)\otimes \delta_{i,j}\; k)
\oplus
\Ext^m_G(T(b-1), \Delta(a)\otimes \delta_{i,\bj}\; k)
\end{multline*}
and we keep on going until we get a tilting module that is
$p$-restricted.
We then use corollary \ref{cor:extD}.

\section{Applying to modules for $\SL_2(k)$ part III}\label{sect:III}

We now want to 
determine $\Ext^q_G(\Delta(pb+j), M\frob \otimes L(i))
\cong \Ext^q_G((M^*)\frob \otimes L(i), \nabla(pb+j))$.

We first explicitly 
construct the $G_1$ injective resolution of $\nabla(pa+i)$.
This is easy using the dual of sequence \eqref{xanth2}
as the middle term is projective as a $G_1$-module.
We have 
$$ 
0\to 
\nabla(pa+i) 
\xrightarrow{\delta_{-1}}
I_0 \xrightarrow{\, \delta_0\, }
I_1 \xrightarrow{\, \delta_1\,}
\cdots  \xrightarrow{\delta_{m-1}}
I_m  \xrightarrow{\delta_{m}}
I_{m+1} \xrightarrow{\delta_{m+1}} \cdots  
$$
with   
$$
I_m \cong \left\{\begin{array}{ll}
                    \nabla(m+b)\frob \otimes Q(j) &\mbox{if $m$ even}\\
                    \nabla(m+b)\frob \otimes Q(\bj) &\mbox{if $m$ odd.}
                 \end{array} \right.
$$
We define $M_m = \ker\delta_{m}$. We have
$$
M_m \cong \left\{\begin{array}{ll}
                    \nabla(p(m+b)+j)   &\mbox{if $m$ even}\\
                    \nabla(p(m+b)+\bj)   &\mbox{if $m$ odd.}
                 \end{array} \right.
$$

A similar argument to that in section
\ref{sect:I} shows that 
$$
\Hom_{G_1}((M^*)\frob \otimes L(i), I_m) 
\cong
\Hom_{G_1}((m^*)\frob \otimes L(i), M_m). 
$$
Thus the induced differential $\delta_{m}^*$ on the $\Hom_{G_1}$
complex is zero.

We have for $m \ge 1$
\begin{align*}
\Ext^m_{G_1}( (M^*)\frob \otimes L(i), \nabla(pb+i)) 
          &\cong \Ext_{G_1}^1((M^*)\frob \otimes L(i), M_{m-1}) \\
&\cong \left\{\begin{array}{ll}
          M\frob\otimes \nabla(m+b)\frob &\mbox{if $m$ even}\\
          0              &\mbox{otherwise}  
                 \end{array} \right.
\end{align*}
\begin{align*}
\Ext^m_{G_1}((M^*)\frob \otimes L(i), \nabla(pb+\bi)) 
          &\cong \Ext_{G_1}^1((M^*)\frob \otimes L(i), M_{m-1}) \\
&\cong \left\{\begin{array}{ll}
          M\frob\otimes \nabla(m+b)\frob &\mbox{if $m$ odd}\\
          0              &\mbox{otherwise}  
                 \end{array} \right.
\end{align*}
and we may thus deduce the following theorem using lemma
\ref{lem:dzero}.

\begin{thm}\label{thm:DL}
Let $M$ be in $\Mod(G)$, $b \in \N$ and $0 \le i \le p-2$.
If $p\ge3$ then 
$$
\Ext^q_{G}(\Delta(pb+i), M\frob\otimes L(i))
\cong \bigoplus_{\substack{n \mbox{\rm{\scriptsize{ even}}}\\ \\{0 \le
      n \le q}}}
      \Ext^{q-n}_G(\Delta(n+b), M)
$$
and
$$
\Ext^q_{G}(\Delta(pb+i), M\frob\otimes L(\bi))
\cong \bigoplus_{\substack{n \mbox{\rm{\scriptsize{ odd}}}\\ \\{0 \le
      n \le q}}}
      \Ext^{q-n}_G(\Delta(n+b), M)
$$
%

If $p=2$ then
$$
 \Ext^q_{G}(\Delta(2b), M\frob) 
\cong
 \bigoplus_{n=0}^{n=q} \Ext^{q-n}_G(\Delta(n+b), M).
$$
\end{thm}
\begin{cor}
Let $b \in \N$, $0 \le i \le p-2$, $q \in \N$ $M \in \Mod(G)$.
If $p\ge3$ then we have
$$
\Ext^q_{G}(\Delta(pb+i), M\frob\otimes L(i)) 
\cong \Ext_{G}^{q-1}(\Delta(p(b+1)+\bi), M\frob\otimes L(i))
\oplus \Ext_{G}^{q}(\Delta(b), M)
$$
$$
\Ext_{G}^{q}(\Delta(pb +i), M\frob\otimes L(\bi))
\cong \Ext_{G}^{q-1}(\Delta(p(b+1)+\bi), M\frob\otimes L(\bi))
$$
where $\Ext^{-1}$ is interpreted as the zero module.
%

If $p=2$ then
$$
 \Ext^q_{G}(\Delta(2b), M\frob) 
\cong
 \Ext^{q-1}_G(\Delta(2b+2), M\frob) \oplus 
 \Ext^q_{G}(\Delta(b), M).
$$
where again, $\Ext^{-1}$ is interpreted as the zero module.
\end{cor}

Thus if $M = \nabla(a)$, $\Delta(a)$ or $L(a)$ we may now completely determine
$\Ext^q_G(\Delta(pb + j),M\frob \otimes L(i))$ and
$\Ext^q_G(M\frob \otimes L(i), \Delta(pb + j))$ using the results of the
previous sections. (This is as
$\Ext^q_G(\Delta(a),M)$ will be zero if $q$ is larger than the highest
weight of $M$ divided by $p$ using lemma \ref{propn:vanish}.)

\section{The quantum case}
All of the results may now be readily generalised to the quantum case,
where we use the quantum group of \cite{dipdonk}.
We need to remember that the Frobenius morphism now goes from the classical
group $\GL_2(k)$ to the quantum group $q$-$\GL_2(k)$ and so when we
`untwist' a twisted module we get a module for $\GL_2$ and not
$q$-$\GL_2(k)$. The $G_1$ cohomology looks slightly different as we
need to keep track of the degree of the modules. Most of the results
of section \ref{sect:sl2g1} have already been generalised to the quantum
case and appear in \cite{cox1} and \cite{coxerd}. The other results may
be proved in a similar fashion. Thus the resolutions we constructed in
sections \ref{sect:I} and \ref{sect:III} may be generalised to the
quantum case and we get essentially the same spectral sequences
appearing. Thus the quantum version of sections \ref{sect:I} to
\ref{sect:III} all hold.
Below are stated the quantum versions of the main corollaries of these
sections for the convenience of the reader.
We will distinguish the modules for classical $\GL_2(k)$ and the
quantum group $q$-$\GL_2(k)$ by putting a bar on the modules for the 
classical groups. We define $\bi = l-i-2$ where $q$ is a primitive
$l$th root of unity.

\begin{thm}
Let $\overline{N} \in\Mod(\GL_2(k))$ be indecomposable and let $b$ be
the highest weight of $\overline{N}$.
For $a \ge b$ with $a-b$ odd, $0 \le i \le l-2$ and $m \in \N$ 
we have
\begin{align*}
\Ext&_{q-\GL_2(k)}^m(\overline{N}\frob\otimes L(\bi+d,d),\Delta(la+i,0))\\
&\cong \Ext_{q-\GL_2(k)}^{m-1}
(\overline{N}\frob\otimes L(\bi+d,d), \Delta(la-1,i+1))
\oplus
\Ext_{\GL_2(k)}^m(\overline{N}\otimes L(f,f), \bDelta(a-1,0))
\end{align*}
where $\Ext^{-1}$ is interpreted as the zero module,
$f=\frac{a-b-1}{2}$ and $d = lf +i+1$.

For $a \ge b$ with $a-b$ even, $0 \le i \le l-2$ and $m \in \N$ we
have
$$
\Ext_{q-\GL_2(k)}^{m}(\overline{N}\frob \otimes L(i+d,d), \Delta(la+i,0))
\cong \Ext_{q-\GL_2(k)}^{m-1}(\overline{N}\frob\otimes L(i+d,d), \Delta(la-1,i+1))
$$
where $d = l(\frac{a-b}{2})$.
\end{thm}

\begin{thm}
For $a \ge b$ with $a-b$ odd, $0 \le i \le l-2$ and $m \in \N$ we
have
\begin{align*}
\Ext&_{q-\GL_2(k)}^m(\Delta(lb+\bi+d,d),\Delta(la+i,0))\\
&\cong \Ext_{q-\GL_2(k)}^{m-1}(\Delta(lb+\bi+d,d), \Delta(la-1,i+1))
\oplus
\Ext_{\GL_2(k)}^m(\bDelta(b+f,f), \bDelta(a-1,0))
\end{align*}
where $\Ext^{-1}$ is interpreted as the zero module,
$f=\frac{a-b-1}{2}$ and $d = lf +i+1$.

For $a \ge b$ with $a-b$ even, $0 \le i \le l-2$ and $m \in \N$ we
have
$$
\Ext_{q-\GL_2(k)}^{m}(\Delta(lb +i+d,d), \Delta(la+i,0))
\cong \Ext_{q-\GL_2(k)}^{m}(\Delta(i+d,d), \Delta(l(a-b)+i,0))
$$
where $d = l(\frac{a-b}{2})$.
If $m\ge 1$ then also 
$$
\Ext_{q-\GL_2(k)}^{m}(\Delta(i+d,d), \Delta(l(a-b)+i,0))
\cong \Ext_{q-\GL_2(k)}^{m-1}(\Delta(i+d,d), \Delta(l(a-b)-1,i+1)).
$$
\end{thm}

\begin{thm}
Let $\overline{M} \in\Mod(\GL_2(k))$ be indecomposable and let $a$ be
the highest weight of $\overline{M}$.
For $a \ge b$ with $a-b$ even, $0 \le i \le l-2$ and $m \in \N$ 
we have
\begin{align*}
\Ext&_{q-\GL_2(k)}^m(\Delta(lb+i+d,d), \overline{M}\frob\otimes L(i,0))\\
&\cong \Ext_{q-\GL_2(k)}^{m-1}
(\Delta(lb+l-1+d,d-l+i+1), \overline{M}\frob\otimes L(i,0))
\oplus
\Ext_{\GL_2(k)}^m(\bDelta(b+f,f), \overline{M})
\end{align*}
where $\Ext^{-1}$ is interpreted as the zero module,
$f=\frac{a-b}{2}$ and $d = lf +i+1$.

For $a \ge b$ with $a-b$ odd, $0 \le i \le l-2$ and $m \in \N$ we
have
\begin{multline*}
\Ext_{q-\GL_2(k)}^{m}(\Delta(lb+i+d,d), \overline{M}\frob \otimes
L(\bi+d,d))\\
\cong \Ext_{q-\GL_2(k)}^{m-1}(\Delta(lb+l-1+d,d-l+i+1), \overline{M}\frob\otimes L(\bi+d,d))
\end{multline*}
where $d = l(\frac{a-b-1}{2})$.
\end{thm}

Block considerations mean that the above cases are the only possible
non-zero $\Ext$ groups of that form.

\section*{Acknowledgements}
I would like to thank Karin Erdmann and the Mathematical Institute in
Oxford for their hospitality in July 2004 
while part of this research was carried out.
Thanks also go to David Benson for some useful discussions
about spectral sequences.
Special thanks also goes to the referee for pointing out some
simplications of the proof of theorems \ref{thm:sectI2} and
 \ref{thm:DL}.


\providecommand{\bysame}{\leavevmode\hbox to3em{\hrulefill}\thinspace}
\providecommand{\MR}{\relax\ifhmode\unskip\space\fi MR }
\providecommand{\MRhref}[2]{%
  \href{http://www.ams.org/mathscinet-getitem?mr=#1}{#2}
}
\providecommand{\href}[2]{#2}

\end{document}